\input amstex
\documentstyle{amsppt}
\magnification=1200
\define \C#1{{\Cal #1}}
\define\po{\parindent 0pt}
\define \p{{\po \it \m Proof. }}
\define \m{\medskip}
\define \ph#1{\phantom {#1}}
\define \T#1{\widetilde {#1}}
\define \Ker{\operatorname{Ker}}
\define \IM{\operatorname {Im}}

\define \cat{\operatorname{cat}}
\define \cwgt{\operatorname{cwgt}}
\define \Crit{\operatorname{Crit}}

\define \Hom{\operatorname{Hom}}
\define \Arn{\operatorname{Arn}}
\define \eps{\varepsilon}
\font\titlefont=cmbx12
\topmatter
\title{\titlefont On the Lusternik--Schnirelmann Category of Symplectic  
Manifolds and the Arnold Conjecture}\endtitle
\rightheadtext{Lusternik--Schnirelmann Category of Symplectic  
Manifolds}
\author{Yuli B. Rudyak and John Oprea}\endauthor
\date June 1997 \enddate
\address{Mathematisches Institut, Universit\"at Heidelberg, 
Im Neuenheimer Feld 288, D-69120 \newline 
\indent Heidelberg 1, Germany}
\endaddress 
\email{july\@mathi.uni-heidelberg.de}\endemail
\address{Department of Mathematics, Cleveland State University,  
Cleveland OH 44115, USA}
\endaddress 
\email{oprea\@math.csuohio.edu}\endemail
\endtopmatter

\head 1. Introduction\endhead

In [Arn, Appendix 9] Arnold proposed a beautiful conjecture  
concerning the relation between the number of fixed points of certain  
(i.e., exact or Hamiltonian) self-diffeomorphisms of a closed  
symplectic manifold $(M,\omega)$ and the minimum number of critical  
points of any smooth (= $C^{\infty}$) function on $M$. The first  
author succeeded in proving this form of the Arnold conjecture [R2]  
under the hypothesis that $\omega$ and $c_1$ vanish on all spherical homology  
classes and that there is equality between the Lusternik--Schnirelmann  
category of $M$ and the dimension of $M$. In this paper, we use a  
fundamental property of category weight to show that, for any closed  
symplectic manifold whose symplectic form vanishes on the image of  
the Hurewicz map, the required equality holds. Thus, we show that the  
original form of the Arnold Conjecture holds for all symplectic  
manifolds having $\omega|_{\pi_2(M)}=0=c_1|_{\pi_2(M)}$.

\head 2. The Arnold Conjecture \endhead

Let $(M^{2n}, \omega)$ be a closed symplectic manifold. A  
symplectomorphism $\phi: M \to M$ (i.e., a diffeomorphism with  
$\phi^*\omega=\omega$) is called {\it Hamiltonian\/} (or {\it exact\/})  
if it belongs to the flow of a time-dependent Hamiltonian vector field on  
$M$. See [HZ] or [MS] for details.

\m 
\noindent{\bf Definition 2.1.} 
\roster
\item"(a)" Define $\Arn(M,\omega)$ to be the  
minimum number of fixed points for any Hamiltonian  
symplectomorphism of $M$.
\item"(b)" Define $\Crit M$ to be the minimum number of critical points  
for any smooth function on $M$.
\endroster

\noindent With these definitions, we can now enunciate the

\m

\noindent{\bf Arnold Conjecture 2.2.} The following inequality 
holds for every closed symplectic manifold $(M, \omega)$: 
$$
\Arn (M,\omega)\geq \Crit M.
$$

\m
The conjecture, usually (but not universally) weakened by replacing  
$\Crit M$ by the cup-length of $M$, has been proved under various  
hypotheses for various classes of manifolds ([CZ], [H], [F1], [F2]).  
The most general approach to the conjecture (in its modified form)  
has been that of Floer who developed a powerful homology theory  
well-suited to the problem. Recently, however, the first author  
proved the following result based on Floer's approach and using  
topological arguments related to a new invariant, the {\it category  
weight\/} (see [R1] for instance).
\proclaim{Theorem 2.3 {\rm ([R2])}} Let $(M,\omega)$ be a closed symplectic manifold, and let $c_1$ be the first Chern class of $M$. Suppose that both $\omega$ and $c_1$ vanish on the image of the  
Hurewicz map 
$
h\: \pi_2(M) \to H_2(M)
$
and that $\cat M =\dim M$.  
Then 
$$
\Arn (M,\omega)\geq \Crit M.
$$
That is, the Arnold Conjecture holds for $M$.
\endproclaim
\m 
Recall that the Lusternik--Schnirelmann category of $M$, $\cat M$,  
is the least integer $n$ such that $M$ may be covered by $n+1$ open  
subsets each of which is contractible inside $M$. Lusternik and  
Schnirelmann ~[LS] showed that
$$
\Crit M \geq \cat M +1
$$
and category has therefore been used as an effective approximation of  
$\Crit M$ ever since. Thus, it is not so surprising that finer  
invariants such as category weight should prove effective in  
analyzing a problem such as the Arnold Conjecture. It is to this new  
invariant that we now turn.

\head 3. Category Weight \endhead

In this section we define {\it category weight\/} and recall some of 
its properties (see [FH], [R1], [S]). Let $H^*(-;G)$ denote singular cohomology with coefficients in an arbitrary abelian group $G$.
\m 
\noindent{\bf Definition 3.1.}  Let $u\in H^*(X;G)$ and denote the {\it  
category weight\/} of $u$ by $\cwgt u$. Then, we say that 
$\cwgt u\geq k$ if $u|_A=0$ for all $A\subseteq X$ with $\cat _X A<k$.
\m
\noindent Here, $\cat_XA$ denotes the relative Lusternik--Schnirelmann  
category defined as the least integer $n$ such that $A$ is covered by  
$n+1$ sets open and contractible in $X$.
\m 
Category weight was introduced recently in [FH] and has been  
developed and applied even more recently in [R1], [R2] and [S]. 
The only properties of $\cwgt$ that we shall need are the following:
\m 
\noindent{\bf Properties of Category Weight 3.2} (see [FH] and 
[R1]).
\roster
\item $\cat X \geq \cwgt x$ for every non-zero $x\in H^*(X;G)$;
\item If $X$ is a $CW$-space then $\cwgt x \leq \deg x$ for every non-zero $x\in H^*(X;G)$;
\item Let $X$ be a metrizable space. Then $\cwgt x \geq 2$ if 
$$
x\in \Ker (\eps^*: H^*(X;G) \to  
H^*(S\Omega X;G)),
$$ 
where $\eps: S \Omega X \to X$ is adjoint to the  
identity on $\Omega X$;
\item $\cwgt(u_1\cdots u_k)\geq \sum_{i=1}^k \cwgt u_i$ for every  
nontrivial cup-product $u_1\cdots u_k$, where 
$u_1, \ldots, u_k\in H^*(X;R)$ and $R$ is any coefficient ring.
\endroster

\head 4. The Main Result \endhead

\proclaim{Theorem 4.1} Let $X$ be a connected $CW$-space, and let  
$u\in H^2(X;\Bbb R)$ be a non-trivial class which vanishes on the  
image of the Hurewicz homomorphism 
$$
h\: \pi_2(X) \to H_2(X).
$$ 
Then  
$\cwgt u=2$.
\endproclaim

\p 
First, recall that, by the Universal Coefficient Theorem, 

$$
H^2(X;\Bbb R) \cong \Hom (H_2(X),\Bbb R).
$$ 

\noindent Hence, the class $u$ may be thought of as a homomorphism from  
$H_2(X)$ to $\Bbb R$. It therefore makes sense to speak about $u$  
vanishing on $\IM(h\: \pi_2(X) \to H_2(X))$ where $h$ is the Hurewicz  
homomorphism. This is usually denoted by $u|_{\pi_2(X)}=0$.
\m By Property 3.2(2) above, we see that $\cwgt u \leq 2$. We shall now  
show that, in fact, equality obtains by demonstrating that $u$  
satisfies the criterion of Property 3.2(3). To do this, recall that  
Hopf's Theorem says that the classifying map
$$
X @>f>> K(\pi_1(X),1)
$$
inducing an isomorphism $f_{\#}\: \pi_1(X) @> \cong >> \pi_1(X)$ also  
gives isomorphisms
$$
f_*\: H_1(X) @> \cong >> H_1(K(\pi,1)) \text{ and } \overline  
f_*\:\frac{H_2(X)}{\IM h} @> \cong >> H_2(K(\pi,1))
$$
where $\pi$ denotes $\pi_1(X)$ and $h\: \pi_2(X) \to H_2(X)$ is the  
Hurewicz homomorphism. (For a straightforward proof, see [LO].) Now,  
$u\in H^2(X;\Bbb R) \cong \Hom (H_2(X);\Bbb R)$ and $u|_{\pi_2(X)}=0$. 
Hence, 
$$
u\in\Hom\left(\frac{H_2(X)}{\IM h},\Bbb R\right) \cong  
\Hom(H_2(K(\pi,1)), \Bbb R) \cong H^2(K(\pi,1);\Bbb R).
$$
We denote the corresponding element of the right-hand group by  
$u_{\pi}$. So, $f^*u_{\pi}=u$.
\m Now consider the commutative diagram
$$
\CD
S\Omega X @>S\Omega f>> S\Omega K(\pi,1)\\
@V\eps_X VV @VV\eps_K V\\
X@>f>> K(\pi,1)@..
\endCD
$$
Clearly, $\Omega K(\pi,1)$ is homotopy equivalent to a discrete  
space, so that $S\Omega K(\pi,1)$ is homotopy equivalent to a wedge  
of circles and, therefore, has no degree 2 cohomology. Hence  
$\eps_K^*=0$ in degree 2. But then
$$
\eps^*_Xu=\eps_X^*(f^*u_{\pi})=(S\Omega f)^*\eps^*_K(u_{\pi})=0.
$$
Thus, by Property 3.2(3), $\cwgt u=2$.
\qed

\proclaim{Corollary 4.2} Let $(M^{2n}, \omega)$ be a symplectic  
manifold with $\omega|_{\pi_2(M)}=0$. Then $\cat M =2n =\dim M$.
\endproclaim

\p By the Theorem, $\cwgt \omega =2$. (Here we write $\omega$ for the  
real cohomology class as well as for the symplectic form). Since  
$(M^{2n}, \omega)$ is a symplectic manifold, the $n$-th wedge product  
of $\omega$ is a volume form and, thus, $\omega^n \neq 0$ in  
cohomology as well. Properties 3.2(1) and 3.2(4) then give
$$
\cat (M) \geq \cwgt \omega^n \geq n \cwgt \omega =2n.
$$
But $\dim M =2n$ and $\dim M \geq \cat M$ in general. Thus, all the  
inequalities are, in fact, equalities and $\cat M =\dim M$.
\qed
\proclaim{Corollary 4.3} Let $(M^{2n}, \omega)$ be a symplectic  
manifold with 
$$
\omega|_{\pi_2(M)}=0=c_1|_{\pi_2(M)}.
$$ 
Then the Arnold Conjecture 
$$
\Arn(M, \omega) \geq \Crit M,
$$
holds for $M$.
\endproclaim

\p Apply Corollary 4.2 to Theorem 2.3.
\qed

\head 5. The Ganea Conjecture. \endhead

In [Ga], Ganea conjectured that the Lusternik--Schnirelmann category  
of the product of a space $X$ with a sphere $S^n, n>0$ should always  
have the category of $X$ plus one. This conjecture has been shown to  
be true under the hypothesis that $\cat X = \dim X$ (see e.g. [R1]).  
Hence, we have 

\proclaim{Theorem 5.1} Let $(M^{2n}, \omega)$ be a symplectic  
manifold with $\omega|_{\pi_2(M)}=0$. Then, for all $n>0$,
$$
\cat (M \times S^n)=\cat M+1.
$$
That is, the Ganea conjecture holds for $M$.
\endproclaim

\head Acknowledgments \endhead

The authors would like to thank the Stefan Banach International  
Mathematics Center in Warsaw for its hospitality and its support of  
the Conference on Homotopy and Geometry, June 9--13, 1997. This paper  
was conceived and written during this conference. The authors would  
also like to thank Greg Lupton for stimulating discussions concerning  
this work.

\newpage

\head References \endhead 

\halign{{\bf #\ }\hfil & \vtop{\parindent0pt
\hsize=31.1em
\hangindent0em\strut#\strut}\cr
{[Ar]} & Arnold, V.I. Mathematical Methods in Classical Mechanics,  
Springer, Berlin Heidelberg New York 1989.
\cr\cr {[CZ]} & Conley, C., Zehnder, E. The Birkhoff--Lewis Fixed Point  
Theorem and a Conjecture of V.I. Arnold, Invent. Math. {\bf 73}  
(1983), 33-49. 
\cr\cr {[FH]} & Fadell, E., Husseini, S. Category weight and Steenrod  
operations, Boletin de la Sociedad Matem\'atica Mexicana, {\bf 37}  
(1992), 151--161.
\cr\cr {[F1]} & Floer, A., Cuplength estimates on Lagrangian  
intersections, Comm. Pure Appl. Math. {\bf 42} (1989), 335--356.
\cr\cr {[F2]} & Floer, A., Symplectic fixed points and holomorphic  
spheres, Commun. Math. Phys. {\bf 120} (1989), 575--611.
\cr\cr {[Ga]} & Ganea, T. Some problems on numerical homotopy  
invariants, Symposium in Algebraic Topology, Seattle 1971, 23--30,  
Lecture Notes in Mathematics {\bf 249}, Springer, Berlin (1971).
\cr\cr {[H]} & Hofer, H. Lusternik--Schnirelmann theory for Lagrangian  
intersections, Annales de l'inst. Henri Poincar\'e-- analyse  
nonlineare, {\bf 5} (1988), 465--499.
\cr {[HZ]} & Hofer, H., Zehnder, E. Symplectic Invariants and  
Hamiltonian Dynamics, Birkh\"auser, Basel, 1994.
\cr\cr {[LO]} & Lupton, G., Oprea, J. Cohomologically Symplectic Spaces:  
Toral Actions and the Gottlieb Group, Trans. Amer. Math. Soc. {\bf  
347}, 1, (1995), 261--288.
\cr\cr {[LS]} & Lusternik, L.A., Schnirelmann, L.G. Methodes  
topologiques dans le probl\`emes variationels, Hermann, Paris (1934).
\cr\cr {[MS]} & McDuff, D., Salamon, D. Introduction to Symplectic  
Topology, Clarendon Press, Oxford 1995.
\cr\cr {[R1]} & Rudyak, Yu. B. On category weight and its applications,  
Preprint 1996.
\cr\cr {[R2]} & Rudyak, Yu. B. On analytical applications of stable  
homotopy (the Arnold conjecture, critical points), Preprint 1997.
\cr\cr {[S]} & Strom, J. Essential category weight, Preprint 1997.\cr}
\enddocument